\documentclass[12pt]{amsart}
\usepackage{amssymb,amsmath}
\usepackage{geometry}    
\usepackage{mathrsfs}

\usepackage{graphicx}

\usepackage{vmargin}
\setmargrb{1in}{1in}{1in}{1in}

\newtheorem{theorem}{Theorem}[section]

\theoremstyle{definition}

\newcommand{\ux}{\underline{\xi}}

\newcommand{\Rn}{\mathbb{R}^n}

\newcommand{\Dn}{\mathbb{D}_n^{\ast} }

\begin{document}

\title[Dirichlet spectrum for one linear form]{Dirichlet spectrum for one linear form}

\author{Johannes Schleischitz}


\thanks{Middle East Technical University, Northern Cyprus Campus, Kalkanli, G\"uzelyurt \\
    johannes@metu.edu.tr ; jschleischitz@outlook.com}


\begin{abstract}
   For $n\geq 2$, we determine the Dirichlet spectrum in $\Rn$ with respect to a linear form and the maximum norm 
   as the entire interval $[0,1]$. This natural result improves
   on recent work of Beresnevich, Guan, Marnat, Ram\'irez and Velani,
   and complements a subsequent paper by the author where
   the analogous result was proved for simultaneous approximation. Various generalisations that can be obtained
   by similar methods as in the
   latter paper are indicated. We believe that our results are an
   important step towards resolving the very open 
   analogous problem for a general system of
   linear forms. 
\end{abstract}

\maketitle

{\footnotesize{

{\em Keywords}: Dirichlet spectrum, linear forms \\
Math Subject Classification 2020: 11J06, 11J13}}



\section{The Dirichlet spectrum for one linear form equals $[0,1]$}


For $\ux=(\xi_1,\ldots,\xi_n)\in\Rn$, let the approximation function
with respect to one linear form be defined as 
\[
\psi_{\ux}^{\ast}(Q)= \min L_{\underline{b}} (\ux)=
 \min | b_0+b_1\xi_1+\cdots+b_n \xi_n |,
\]
where the minimum is taken over all integer vectors $\underline{b}=(b_0,\ldots,b_n)\neq \underline{0}$ with $\max_{1\leq i\leq n} |b_i|\leq Q$. Derive the Dirichlet constant of $\ux$ via
\[
\Theta^{\ast}(\ux)= \limsup_{Q\to\infty} \; Q^n \psi_{\ux}^{\ast}(Q).
\]
Then $\Theta^{\ast}(\ux)\in [0,1]$ by a well-known variant of
Dirichlet's Theorem.
Define
\[
\mathbb{D}_n^{\ast}= \{ \Theta^{\ast}(\ux)\!:\; \ux\in\Rn  \} \subseteq [0,1]
\]
the Dirichlet spectrum with respect to approximation of one linear form 
and the maximum norm. We remark that we 
use the notation with asterisks above to distinguish our setting
from the case of simultaneous approximation from~\cite{am2, ich}.
Our main result is

\begin{theorem} \label{thm1}
	For any $n\ge 2$ we have $\Dn=[0,1]$. More precisely, 
	for any $c\in [0,1]$ there exist uncountably many $\ux\in\Rn$
	with $\Theta^{\ast}(\ux)=c$. Moreover, we may choose them Liouville vectors, i.e. satisfying
	$\liminf_{Q\to\infty} Q^N \psi^{\ast}_{\ux}(Q)=0$ for all $N$.
\end{theorem}

In fact the constructive proof also yields $\ux$ that are Liouville
vectors with respect to simultaneous approximation, 
i.e. with respect to $\psi_{\ux}$ from~\cite{ich} 
in place of $\psi^{\ast}_{\ux}$,
a stronger property. In particular all its coordinates are Liouville numbers.
The claims of the theorem fail for $n=1$, for this case 
we refer to~\cite{am2} for further references.

Theorem~\ref{thm1} refines work of~\cite{beretc} and~\cite[Appendix]{ich}. 
Both papers yield a countable partition of $[0,1]$ 
into intervals each having non-empty intersection with $\Dn$, and
admit (but do not force) the corresponding $\ux\in\Rn$ to be Liouville vectors.
The main purpose of~\cite{beretc} was to show the existence of $\ux$ 
satisfying the considerably weaker claims that $\ux$
has Dirichlet constant in $(0,1)$, i.e. is Dirichlet improvable but not singular, and also not badly approximable, i.e. $\liminf_{Q\to\infty} Q^n\psi_{\ux}^{\ast}(Q)=0$.
It employed Roy's unconstructive existence result on parametric geometry of numbers~\cite{roy} as the essential
tool, which led to aforementioned induced intervals of the 
form $[\kappa_n c,c]$ 
for explicitly computable $\kappa_n\in (0,1)$, which however decays
hyperexponentially with $n$.
In the latter paper~\cite{ich}, the aforementioned partition
was a consequence of~\cite[Theorem~2.1]{ich}  where it was 
constructively shown that
the dual Dirichlet spectrum $\mathbb{D}_m$ with respect to simultaneous approximation to $\xi_1,\ldots, \xi_m$ equals the entire interval
$[0,1]$ for $m\ge 2$, as 
in Theorem~\ref{thm1}, together with a transference result
due to German~\cite{german}. See also~\cite{am2, as, kr} 
on the Dirichlet spectrum with respect to simultaneous approximation.
In contrast to~\cite{beretc}, in
the induced intervals in~\cite{ich}, the quotient of interval
endpoints depends explicitly on the larger endpoint $c$, and this ratio
$\kappa_n(c)$ tends to $0$ like $c^{n^2}$ as $c\to 0$. A comparison between the two results, distinguishing between small and large values of $c$,
was given in~\cite[Appendix]{ich}. We further want to 
mention a related paper of Marnat~\cite{marnat}, which however
does not improve on $\kappa_n$ from~\cite{beretc}.
Theorem~\ref{thm1} finally gives the complete, natural answer  
on the Dirichlet spectrum for one linear form
and thereby closes the gap arising in~\cite{beretc, marnat, ich}.
We note that the problem to determine the Dirichlet spectrum
for a general system of $m$ linear forms in $n$ variables, 
both of which $\Dn$ and $\mathbb{D}_m$
are special cases of, remains open, as does the more moderate
problem formulated in~\cite[Problem~4.1]{beretc} on the existence of Dirichlet
improvable but not badly approximable matrices (even for equal weights). 
We comment a little more on 
this topic in~\S~\ref{co}. 

We can generalize Theorem~\ref{thm1} in various directions, with similar
methods as in~\cite{ich}. This includes exact uniform approximation with respect to
more general functions $\Phi^{\ast}(Q)$ in place of $cQ^{-n}$, metrical
results for packing and Hausdorff dimension for $\ux$ (for Hausdorff dimension we
omit the property of $\ux$ being Liouville), 
generalizations to certain Cantor sets, assertions on the spectrum
with respect to other norms (applied to $\underline{b}$), and additional
refined information on ordinary approximation. We only state
one examplary result.

\begin{theorem}  \label{2}
	Let $n\geq 2$ be an integer.
	Assume $\Phi^{\ast}: \mathbb{N}\to (0,\infty)$ satisfies $\Phi^{\ast}(t)<t^{-n}$
	for $t\ge t_0$ and 
	\begin{equation}  \label{eq:CC}
	\liminf_{\alpha\to 1^{+}}\liminf_{t\to\infty}\frac{\Phi^{\ast}(\alpha t)}{ \Phi^{\ast}(t) } \geq 1.
	\end{equation}
	Then there exist uncountably many $\ux\in\Rn$ satisfying
	\[
	\limsup_{Q\to\infty} \frac{\psi_{\ux}^{\ast}(Q)}{\Phi^{\ast}(Q)}=1.
	\]
\end{theorem}


The claim can be interpreted complementary to the special case of one linear form of a result by Jarn\'ik~\cite{jarnik}.
He constructed vectors $\ux$ with the property $\psi_{\ux}^{\ast}(Q)<\Phi^{\ast}(Q)$ for all large $Q$, for any given $\Phi^{\ast}$ of decay $\Phi^{\ast}(t)<t^{-n}$, but provides no control for the lower bound
of $\psi_{\ux}^{\ast}$. In~\cite{jarnik} algebraic independence of the coordinates of $\ux$ was imposed as well, however by some variational
argument assuming this does not cause any problem in Theorem~\ref{2} either. So only the additional condition \eqref{eq:CC} is a disadvantage of Theorem~\ref{2}, which excludes functions of (local)
exponential decay. The ``very singular'' vectors induced
by very fast decreasing $\Phi^{\ast}$ have frequently occurred 
in the work of Moshchevitin, particularly regarding geometric
aspects of their induced best approximating integer linear forms, we only want to
quote~\cite{mosh} exemplarily here. See also~\cite{ichneu} for
inhomogeneous approximation.

\section{Proof of Theorem~\ref{thm1} }

\subsection{Outline}
Given $c\in[0,1]$, we construct $\ux\in\mathbb{R}^n$ with $\Theta^{\ast}(\ux)=c$.
We may assume strict inequalities $0<c<1$ as it is well-known that the set of singular vectors defined via $\Theta^{\ast}(\ux)=0$ is rather large (see~\cite{dfsu1},~\cite{dfsu2} for example), and similarly 
for $\Theta^{\ast}(\ux)=1$ the according set of $\ux$ has full 
$n$-dimensional Lebesgue measure due to work of Davenport
and Schmidt~\cite{davs}. 
The construction described in~\S~\ref{2.1} is very similar to~\cite[\S~5.1]{ich}
where simultaneous approximation was studied, with only a small twist.
The verification in~\S~\ref{s2.2},~\ref{2.3} below is based on similar concepts
as well, however some new ideas are required, in particular the 
most challenging Case 3 of~\S~\ref{s2.2} employs a fundamentally new technique.

\subsection{Construction of $\ux$ }  \label{2.1}

Define an increasing sequence of positive integers
recursively as follows: 
Define the initial $n$ terms by $a_j=8j!$ 
for $1\leq j\leq n$.
For $k\geq 1$, having constructed the first $nk$ terms $a_1,\ldots,a_{nk}$, let
the next $n$ terms be given by the recursion
\begin{equation}  \label{eq:h}
a_{nk+1}= a_{nk}^{M_k}
\end{equation}
and
\begin{equation}  \label{eq:tats}
a_{nk+2}= a_{nk+1}^2, \quad a_{nk+3}=a_{nk+1}^3, \quad \ldots, \quad a_{nk+n-1}=a_{nk+1}^{n-1},
\end{equation}
and finally
\begin{equation} \label{eq:itus}
a_{nk+n}=a_{n(k+1)}= \lceil c^{-1}a_{nk+1}\rceil\cdot a_{nk+n-1}\in (a_{nk+1}^n,\infty),
\end{equation}
with integers $M_k\to\infty$ that tend to infinity fast enough.

We observe that
\begin{equation} \label{eq:divid}
a_i|a_{i+1}, \qquad i\geq 1.
\end{equation}
Define the components $\xi_j$
of $\ux$ via
\begin{equation}  \label{eq:sinus}
\xi_j= \sum_{h=0}^{\infty} \frac{1}{a_{nh+j}}, \qquad 1\leq j\leq n,
\end{equation}
that is we sum the reciprocals
over the indices congruent to $j$ modulo $n$.
We show that $\Theta^{\ast}(\ux)=c$ for $\ux$ constructed above,
and that they are Liouville vectors. Due to the freedom in choice of the $M_k$ we clearly obtain uncountably many $\ux$.

For the sequel,
for $1\leq j\leq n$ and $k\geq 0$ integers, fix the notation $\xi_j=S_{j,k}+R_{j,k}$ with 
\[
S_{j,k}= \sum_{h=0}^{k} a_{nh+j}^{-1}\in\mathbb{Q},\qquad
R_{j,k}=\sum_{h=k+1}^{\infty} a_{nh+j}^{-1}.
\]
Let $S_{0,k}=1$ and $R_{0,k}=0$ for $k\geq 0$. Then for any 
$\underline{b}\in \mathbb{Z}^{n+1}$
and any $k\geq 1$ we have
\begin{equation} \label{eq:jo}
L_{\underline{b}}(\ux)= |\sum_{j=0}^{n} b_j S_{j,k} + \sum_{j=0}^{n} b_j R_{j,k}|. 
\end{equation}

\subsection{
Proof of $\Theta^{\ast}(\ux)\leq c$ and Liouville property.} \label{s2.2}
Let $Q$ be any large parameter. We need to construct $\underline{b}\in\mathbb{Z}^{n+1}\setminus \{ \underline{0} \}$ 
with $\max_{1\leq j\leq n} |b_j|\leq Q$ and
inducing $L_{\underline{b}}(\ux)<c(1+o(1))Q^{-n}$ as $Q\to\infty$.
There is a unique large integer $k$ such that $a_{nk}\leq Q< a_{n(k+1)}$.
Observing that 
\[
a_{nk}< \frac{ a_{nk+1}a_{n(k-1)} }{ a_{nk} } < a_{nk+1} <a_{n(k+1)} 
\]
we distinguish the three cases for $Q$ given by
\[
Q\in [a_{nk},a_{nk+1}a_{n(k-1)}/a_{nk}),\quad  
Q\in [a_{nk+1}a_{n(k-1)}/a_{nk}, a_{nk+1}),\quad Q\in [a_{nk+1}, a_{n(k+1)}),.
\]
We start with the easiest case and keep the most challenging middle interval for the end.

Case 1: We have $Q\in [a_{nk+1}, a_{n(k+1)})$. Then we may simply take
\[
b_1=a_{nk+1}, \qquad b_2=b_3=\cdots=b_{n}=0,
\]
and $b_0=-b_1S_{1,k}=- \sum_{j\leq k} a_{nk+1}/a_{nj+1}$ the nearest integer to $- b_1\xi_1$, which
by \eqref{eq:divid} is indeed an integer. 
Clearly $\max_{1\leq j\leq n} |b_j|\leq Q$.
From \eqref{eq:jo}, \eqref{eq:h}
and \eqref{eq:itus} we get 
\begin{align*}
L_{\underline{b}}(\ux)&= |b_1S_{1,k}+b_1R_{1,k}+b_0|= |b_{1}R_{1,k}|=
a_{nk+1}\sum_{j\geq k+1} a_{nj+1}^{-1}  \\ &\leq 
2 a_{nk+1} a_{n(k+1)+1}^{-1}\leq 2a_{n(k+1)}^{1/n}\cdot(c^{-1}a_{n(k+1)}^{M_{k+1}})^{-1}
= 2ca_{n(k+1)}^{1/n-M_{k+1} } <ca_{n(k+1)}^{-n}\leq  cQ^{-n}
\end{align*}
as soon as we choose $M_{k+1}$ large enough in the next step 
of the construction. Thus $\psi^{\ast}_{\ux}(Q)\leq cQ^{-n}$ for those $Q$.
If we let $M_k\to\infty$, by taking $Q=a_{nk+1}$ and the $b_j$ as above, it is further easy to see that the induced $\ux$ becomes a Liouville vector.

Case 2: We have $Q\in [a_{nk},a_{nk+1}a_{n(k-1)}/a_{nk})$. 
Then we choose
\[
b_1=b_2=\cdots=b_{n-1}=0, \qquad b_n=a_{nk}
\]
and finally $b_0=- b_nS_{n,k}=- \sum_{j\leq k} a_{nk}/a_{nj}$ the nearest integer to $-b_n\xi_n$, which again
by \eqref{eq:divid} is indeed an integer. It is again obvious that 
$\max_{1\leq j\leq n} |b_j|\leq Q$.
From \eqref{eq:jo} we get
\[
L_{\underline{b}}(\ux)= |b_{n}S_{n,k}+ b_{n}R_{n,k}+b_0|=|b_n R_{n,k}|=
a_{nk}R_{n,k}= \sum_{j\geq k+1} \frac{ a_{nk} }{ a_{nj} }
\leq 2 \frac{ a_{nk} }{ a_{n(k+1)} }. 
\]
We want this to be $< cQ^{-n}$. By assumption on $Q$ it suffices to have
\[
2\frac{ a_{nk} }{ a_{n(k+1)} } < c \left( \frac{ a_{nk+1}a_{n(k-1)} }{ a_{nk} } \right)^{-n}.
\]
Since $a_{n(k+1)}\geq c^{-1}a_{nk+1}^{n}$
by \eqref{eq:itus}, a sufficient condition is
\begin{equation}  \label{eq:jodad}
2a_{n(k-1)}< a_{nk}^{1-1/n}.
\end{equation}
Since $a_{nk}>c^{-1}a_{n(k-1)}^{M_{k-1}n}>a_{n(k-1)}^{M_{k-1}n}$ 
by \eqref{eq:h} and \eqref{eq:itus}, and as $n\geq 2$, this can clearly be achieved for all $k$
by choosing the $M_{j}\geq 2$. Again
we infer $\psi^{\ast}_{\ux}(Q)\leq cQ^{-n}$ for these values of $Q$.

Case 3: Assume finally $Q\in [a_{nk+1}a_{n(k-1)}/a_{nk}, a_{nk+1})$. 
For $0\leq \ell \leq k-1$, define
\[
N_{\ell}=N_{\ell}(k):= a_{nk+1}a_{n\ell}\left( \frac{ 1 }{ a_{n(\ell+1)} }  + \frac{ 1 }{ a_{n(\ell+2)} }
+ \cdots + \frac{ 1 }{ a_{nk} }  \right),
\] 
where we let $a_0=1$. Note that all $N_{\ell}$ are integers by \eqref{eq:divid} and we can assume
that $N_0>N_1> N_2> \cdots>N_{k-1}$, it suffices to choose $M_j\ge 2$. The assumption of Case 3 
translates into $Q\ge N_{k-1}$, hence there exists a minimum
integer $0\leq e\leq k-1$ such that
\[
Q\ge N_{e}.
\]
Let us define our linear form via
\[
b_1=N_e, \qquad b_2=b_3=\cdots=b_{n-1}=0,\qquad b_n=-a_{ne}
\]
and finally put $b_0$ the closest
integer to $-(b_1\xi_1+b_n\xi_n)$. Clearly 
$|b_n|= a_{ne}\leq N_e\le Q$ as well, so $\max_{1\leq j\leq n} |b_j|\leq Q$. 
We may first rewrite $L_{\underline{b}}(\ux)$ from \eqref{eq:jo}
with our present $k$ as
\begin{align}  \label{eq:11}
L_{\underline{b}}(\ux)&=  |\sum_{j=0}^{n} b_j S_{j,k} + \sum_{j=0}^{n} b_j R_{j,k}| \\ \nonumber
&= | b_0+b_1S_{1,k}+ b_nS_{n,k}+ b_1R_{1,k}+ b_nR_{n,k} |  \\ \nonumber
&=| b_0+ b_1S_{1,k-1} + (b_1(S_{1,k}-S_{1,k-1}) + b_nS_{n,k}) + (b_1R_{1,k}+b_nR_{n,k}) |.
\end{align}

Next we observe that $b_1 S_{1,k-1}\in\mathbb{Z}$. Indeed,
we may write it as
\[
b_1 S_{1,k-1}=U_ka_{ne} \cdot \frac{a_{nk+1}}{a_{n(k-1)+1}a_{nk}}, \qquad 
U_k:=\frac{a_{nk}}{a_{n(e+1)}}+ \frac{a_{nk}}{a_{n(e+2)}}+\cdots+\frac{a_{nk}}{a_{nk}}\in\mathbb{Z},
\]
where $U_k\in\mathbb{Z}$ by \eqref{eq:divid}.
But now
$(a_{n(k-1)+1}a_{nk})| a_{nk+1}$ again since by \eqref{eq:divid} 
we have $a_{n(k-1)+1}| a_{nk}$
and by \eqref{eq:itus} we have $a_{nk}^2| a_{nk+1}$,
as we may assume $M_k\ge 2$.
Moreover we easily check via \eqref{eq:divid} that
\[
b_1(S_{1,k}-S_{1,k-1}) + b_nS_{n,k}= N_{e}/a_{nk+1}- a_{ne}S_{n,k}= -a_{ne}(a_n^{-1}+a_{2n}^{-1}+\cdots+a_{ne}^{-1})\in\mathbb{Z}
\]
as well, where we mean that the sum vanishes if $e=0$. Since the error terms $b_1R_{1,k}$ and $b_nR_{n,k}$
are easily seen to be each smaller than $1/4$ by absolute value, we get
\[
b_0=-b_1 S_{1,k-1}-(b_1(S_{1,k}-S_{1,k-1}) + b_nS_{n,k})\in\mathbb{Z}
\]
and can simplify \eqref{eq:11} to
\begin{align}  \label{eq:1}
L_{\underline{b}}(\ux)&=  |b_1R_{1,k}+b_nR_{nk}| \\ \nonumber
\nonumber &= 
|N_eR_{1,k}-a_{ne}R_{n,k}|\\ &=
 a_{ne}\cdot \left| a_{nk+1}\left(\frac{1}{a_{n(e+1)}}+\frac{1}{a_{n(e+2)}}+\cdots+\frac{1}{a_{nk}} \right)\sum_{j\geq k+1} \frac{1}{ a_{nj+1} } -
\sum_{j\geq k+1} \frac{1}{ a_{nj} }\right|. \nonumber
\end{align}
We distinguish two cases again.

Subcase 3i: Assume $e\geq 1$. Choosing $M_{k+1}\geq 2$, 
from \eqref{eq:h} we easily
see that $a_{nk+1}/a_{n(k+1)+1}=o(a_{n(k+1)}^{-1})$ as $k\to\infty$
and consequently $a_{ne}/a_{n(k+1)}$ is the dominating term 
in \eqref{eq:1}. 
It suffices to estimate \eqref{eq:1}
by
\[
L_{\underline{b}}(\ux)\leq
2\frac{a_{ne}}{ a_{n(k+1)} }.
\]
Again  we want this to be $<cQ^{-n}$. 
This is equivalent to
\[
Q< \left( \frac{ca_{n(k+1)}}{2a_{ne}}  \right)^{1/n}.
\]
By assumption $e\ge 1$ and minimality of $e$, we have
\[
Q\leq N_{e-1}< 2\frac{ a_{n(e-1)} }{ a_{ne} } a_{nk+1}.
\]
Hence, it suffices to confirm
\[
2\frac{ a_{n(e-1)} }{ a_{ne} } a_{nk+1}\leq \left( \frac{ca_{n(k+1)}}{2a_{ne}}  \right)^{1/n}.
\] 
But using $a_{n(k+1)}\geq c^{-1}a_{nk+1}^{n}$ from \eqref{eq:itus} and $n\geq 2$,
a sufficient condition for the latter is
\begin{equation}  \label{eq:indeed}
a_{n(e-1)}\leq \frac{1}{4} \cdot a_{ne}^{1/2}\leq \frac{1}{4} \cdot a_{ne}^{1-1/n}.
\end{equation}
Now we conclude similar as for \eqref{eq:jodad}:
For $e=1$ we may assume this
since $a_0=1$ and $a_n=8n!\geq 16$, 
and for $e>1$ by \eqref{eq:tats} and \eqref{eq:h} we have
\[
a_{ne}> c^{-1}a_{n(e-1)}^{M_{e-1}n}> a_{n(e-1)}^{M_{e-1}n}\geq a_{n(e-1)}^{2M_{e-1} }.
\]
If we let $M_j\geq 2$ for all $j\geq 1$
and choose larger initial terms $M_j$ if necessary, 
the quotient $a_{n(e-1)}/a_{ne}^{1/2}$ from
\eqref{eq:indeed} can be made arbitrarily small, uniformly in 
$c\in [0,1]$ and $e\geq 1$.
Thus again $\psi^{\ast}_{\ux}(Q)\leq cQ^{-n}$ for $Q$ as above
in any case.

Subcase 3ii: Finally assume $e=0$, i.e. $Q>a_{nk+1}(a_{n}^{-1}+a_{2n}^{-1}+\cdots+a_{kn}^{-1})$. Again
$a_0/a_{n(k+1)}=a_{n(k+1)}^{-1}$ is the dominating term 
in \eqref{eq:1}, for large $k$. 
Thus, from \eqref{eq:itus} as $Q\to\infty$ (or equivalently $k\to\infty$) we get
\[
L_{\underline{b}}(\ux)\leq a_{n(k+1)}^{-1}(1+o(1))=ca_{nk+1}^{-n}(1+o(1))\leq cQ^{-n}(1+o(1)).
\]
This yields $\psi^{\ast}_{\ux}(Q)\leq c(1+o(1))Q^{-n}$ as $Q\to\infty$ 
for values $Q$ in question.
Since we exhausted all cases,
we infer $\Theta^{\ast}(\ux)\leq c$.

\subsection{
Proof of $\Theta^{\ast}(\ux)\geq c$.} \label{2.3}
Take $Q=a_{nk+1}-1$ for large $k$. Consider any non-trivial linear form $L_{\underline{b}}$ 
of height at most $Q$.
Let $i\in \{ 1,2,\ldots, n\}$ be the largest index with $b_i\neq 0$.
By \eqref{eq:divid}, the denominator of any $S_{j,k}\in\mathbb{Q}$ written in lowest terms divides $a_{nk+j}$. Assume for the moment that
for $j=i$, the denominator of $S_{i,k}$ equals $a_{nk+i}$ when reduced.
If $i=1$, then $b_1S_{1,k}\notin\mathbb{Z}$ by $|b_1|\leq Q<a_{nk+1}$, thus also 
$b_0+b_1S_{1,k}\notin\mathbb{Z}$.
If $i>1$,
then since $|b_i|\leq Q<a_{nk+1}\leq a_{nk+i}/a_{nk+i-1}$ 
by \eqref{eq:tats}, \eqref{eq:itus}
(with equality in the last inequality
if $i<n$), we have
that $b_i S_{i,k}$ in reduced form still has denominator strictly larger than $a_{nk+i-1}=(i-1)a_{nk+1}$.
On the other hand, by \eqref{eq:divid} all other $b_jS_{j,k}$ for 
$0\leq j<i$ have denominator dividing $(i-1)a_{nk+1}$, after reduction.
Hence this also applies to their sum. 
In both cases $i=1$ or $i>1$, we infer that the entire sum 
\[
\sum_{j=0}^{n} b_jS_{j,k} =  \sum_{j=0}^{i} b_jS_{j,k} = b_i S_{i,k}+
\sum_{j=0}^{i-1} b_jS_{j,k} 
\] 
cannot vanish. But again by \eqref{eq:divid} it can be written
as a rational with denominator dividing $a_{nk+i}$ when reduced.
Therefore its modulus is bounded from below by 
\[
\left|\sum_{j=0}^{n} b_jS_{j,k}\right|\geq a_{nk+i}^{-1}\ge a_{n(k+1)}^{-1}> c(1-o(1))a_{nk+1}^{-n}> c(1-o(1))Q^{-n},
\]
as $Q\to\infty$ (or equivalently $k\to\infty$), where we
used \eqref{eq:itus} and that $Q$ is very close to $a_{nk+1}$. The sum of remainder terms can be estimated via
\[
\left|\sum_{j=0}^{n} b_jR_{j,k}\right|\leq nQ \max_{1\leq j\leq n} |R_{j,k}|\ll_{n} a_{nk+1}a_{n(k+1)+1}^{-1}=o(a_{nk+1}^{-n})=o(Q^{-n}),
\] 
as $Q\to\infty$ as soon as we let all $M_j\geq 2$ by \eqref{eq:h}. 
By triangle inequality and \eqref{eq:jo} we get
\begin{align*}
L_{\underline{b}}(\ux)&=
\left| \sum_{j=0}^{n} b_jS_{j,k} + \sum_{j=0}^{n} b_jR_{j,k}\right| 
\\ &\geq \left| \sum_{j=0}^{n} b_jS_{j,k}\right| - \left| \sum_{j=0}^{n} b_jR_{j,k}\right|\geq
c(1-o(1))Q^{-n}, \qquad Q\to\infty.
\end{align*}
Thus, if we let $k\to\infty$, we get $\Theta^{\ast}(\ux)\geq c$ as desired.

It remains to be shown that $S_{i,k}$ is already reduced with
denominator $a_{nk+i}$, for $1\leq i\leq n$. We use induction on $k$. 
It is clear that $S_{i,0}= 1/a_{i}$ is reduced.
Now assume $S_{i,k-1}=u/v$ is in lowest
terms with denominator $v=a_{n(k-1)+i}$. If $i<n$ 
by \eqref{eq:h}, \eqref{eq:tats} and \eqref{eq:itus} we have that
\[
S_{i,k}= \frac{u}{v}+ \frac{1}{ a_{nk+i} } = \frac{u}{v}+ \frac{1}{(vG_{i,k})^{iM_k }} = \frac{  v^{iM_k-1}G_{i,k}^{iM_k}u +1  }{ v^{iM_k}G_{i,k}^{iM_k} }= \frac{  v^{iM_k-1}G_{i,k}^{iM_k}u +1  }{ a_{nk+i} }
\]
for some integer $G_{i,k}=a_{n(k-1)+1}^{n-1-i}\lceil c^{-1}a_{n(k-1)+1}\rceil$.
This is clearly again reduced with denominator $a_{nk+i}$ since upon $M_k\geq 2$ the numerator is congruent to $1$ modulo any prime
divisor of $vG_{i,k}$, the latter being an integer multiple of the radical of the denominator. 
Finally in case of $i=n$, we have
\[
S_{n,k}= \frac{u}{v}+ \frac{1}{ a_{nk+n} } = \frac{u}{v}+ \frac{1}{H_{n,k}v^{M_k (n-1)}} = \frac{  v^{M_k (n-1)-1}H_{n,k}u +1  }{ v^{M_k (n-1)}H_{n,k}}= \frac{  v^{M_k (n-1)-1}H_{n,k}u +1  }{ a_{nk+n} }
\]
where $H_{n,k}=\lceil c^{-1}a_{nk+1} \rceil=\lceil c^{-1}v^{M_{k}}\rceil$, and a similar argument applies.

\section{Sketch of the proof of Theorem~\ref{2} and concluding remarks} \label{co}

To prove Theorem~\ref{2}, we replace \eqref{eq:itus} by
\[
a_{nk+n}= a_{nk+n-1}\cdot L_k^{\ast}, \qquad k\geq 1,
\]
with 
\[
L_k^{\ast}=  \min \{ z\in \mathbb{N}:\; (za_{nk+n-1})^{-1}<\Phi^{\ast}(a_{nk+1}) \}.
\]
The adaptions in the proof above can be readily made, very similarly as in~\cite[\S~5.2]{ich}. We omit the details.

The proof of $\Theta^{\ast}(\ux)\leq c$ above shows that independent of $n$ and $Q$, we only require at most
$3$ non-vanishing coefficients $b_j$ in the integer vector $\underline{b}$ 
to guarantee the estimate $\psi^{\ast}_{\ux}(Q)<(c+o(1))Q^{-n}$, as $Q\to\infty$.
However, it is probably not true that these sparse linear forms are
minimal points to $\ux$ in general (i.e. realize the minimum $\psi_{\ux}^{\ast}(Q)$ over all $\underline{b}$ with norm in the according intervals).
The proof shows that upon taking larger initial terms $a_1,\ldots,a_n$ if necessary, actually taking integers $M_k\geq 2$ in each step 
in \eqref{eq:h} should be sufficient. Then our vector $\ux$ does not consist of Liouville number entries, in fact it is no Liouville vector in the sense of Theorem~\ref{thm1}.

It would be highly desirable to combine ideas from~\cite{ich}
and the present note to generalize the construction
to a system of $m$ linear forms in $n$ variables. 
However, particularly the modification of our integer vectors
in Case 3 to the case $m>1$ causes technical problems that we
have not yet been able to overcome.

\vspace{0.5cm}

{\em The author thanks the referee for the careful reading and suggestions.}

\end{document}